\documentclass{elsarticle}

\usepackage{amsmath}
\usepackage{amsfonts} 
\usepackage{amsbsy} 
\usepackage{amssymb}
\usepackage{bm}
\usepackage{accents}

\newtheorem{thm}{Theorem}

\newproof{pf}{Proof}

\newcommand{\const}{\mathop{\rm const}\nolimits}

\journal{arXiv} 

\begin{document}

\begin{frontmatter}

\title{Some methods for solving equations with an operator function and applications for problems with a fractional power of an operator}

\author[nsi,uni]{Petr N. Vabishchevich\corref{cor}}
\ead{vabishchevich@gmail.com}

\address[nsi]{Nuclear Safety Institute, Russian Academy of Sciences, Moscow, Russia}
\address[uni]{North-Eastern Federal University, Yakutsk, Russia}

\cortext[cor]{Corresponding author}

\begin{abstract}

Several applied problems are characterized by the need to numerically solve equations with an operator function (matrix function).
In particular, in the last decade, mathematical models with a fractional power of an elliptic operator and numerical methods for their study have been actively discussed.
Computational algorithms for such non-standard problems are based on approximations by the operator function.
The most widespread are the approaches using various options for rational approximation.
Also, we note the methods that relate to approximation by exponential sums.
In this paper, the possibility of using approximation by exponential products is noted.
The solution of an equation with an operator function is based on the transition to standard stationary or evolutionary problems.
General approaches are illustrated by a problem with a fractional power of the operator.
The first class of methods is based on the integral representation of the operator function under rational approximation, approximation by exponential sums, and approximation by exponential products.
The second class of methods is associated with solving an auxiliary Cauchy problem for some evolutionary equation.

\end{abstract}


\end{frontmatter}

\section{Introduction}

Classical applied mathematical models \cite{lin1988mathematics}  are based on systems of elliptic, parabolic and hyperbolic equations.
Space dependence is local in nature and is described by elliptic operators \cite{evans2010PDE}.
Recently, increased attention has been paid to nonlocal mathematical models that are associated with fractional derivatives \cite{bisci2016variational}.
In particular, nonlocal diffusion processes are modeled by fractional powers of elliptic operators \cite{Pozrikidis2016}. 

In general, we can consider problems with an operator function.
At the discrete level, we get problems for function of a matrix \cite{higham2008functions}; in this case we are looking for a solution $u$ of the equation $f(A) u = \varphi$.
In applied problems, the matrix $A$ is associated with some discretization of the elliptic operator.
With the standard finite-difference or finite-element approximation, we have a sparse matrix $A$.
The main computational problem is related to the fact that $f(A)$ is a full matrix with hard-to-compute elements.

Various methods are being developed to find an approximate solution to  $f^{-1}(A) \varphi$. 
The main approach is based on some approximation of the function $f^{-1}(\lambda)$: $r(\lambda) \approx f^{-1}(\lambda)$. 
In this case, the approximate solution is $r(A) \varphi$.
The choice of the approximation $r(A)$ is subject to the requirement of the computational implementation acceptability.
Ideally, calculating $r(A) \varphi$ should not be more difficult in principle than calculating $A^{-1} \varphi$; the computational complexity of solving the problem $f(A) u = \varphi$ must be comparable to the solution of the standard problem $A u = \varphi$.

When approximating functions, rational approximations are most widely used \cite{braess1986nonlinear}. 
This approach is actively used in various versions for the approximate solution of problems with a fractional power of an elliptic operator (see, for example, \cite{harizanov2020rev} and the cited literature).
We can highlight the approximation by exponential sums \cite{braess1986nonlinear} separately.
Some results on using this approach for power function approximation are presented in \cite{beylkin2010approximation}.

This paper considers some rational approximations and approximations by exponential sums when solving equations with an operator function for self-adjoint positive definite operators in a finite-dimensional Hilbert space.
We will touch  only a few aspects of the theory.
The accuracy of the approximate solution is discussed for the absolute and relative accuracy of the approximation of functions.
The influence of the error in solving auxiliary problems in the computational implementation of the operator function approximation procedures is investigated separately.
We have identified a new class of approximations using the product of exponentials.

The mentioned approaches for the approximate solution of equations with an operator function are illustrated for problems with fractional operator powers.
The focus is on constructing rational approximations, approximations by exponential sums, and approximations by exponential products.
In the literature, the integral representation of the operator function is most commonly used.
In this case, we can associate specific approximations of the operator function with the corresponding quadrature formulas.
The possibilities of constructing an approximate solution based on a particular Cauchy problem for evolution equations are also noted.

The paper is organized as follows.
In Section 2, we set up notation and terminology. We formulate an approximate solution to an equation with an operator function for self-adjoint positive definite operators in a finite-dimensional Hilbert space.
In Section 3, our main results are stated and proved.
We investigate the accuracy of the approximate solution using rational approximation, approximation by exponential sums, and approximation by exponential products.
In Section 4, we use integral representations for fractional operator approximations.
Section 5 is devoted to the Cauchy problem methods for auxiliary evolutionary equations.
Concluding remarks are drawn in Section 6.

\section{Formulation of the problem} 

Let $H$ be a finite-dimensional Hilbert space. 
The scalar product for $u, v \in H$ is $(u, v)$, and the norm is $\| u \| = (u, u)^{1/2}$.
For a self-adjoint and positive definite operator $D$, we define a Hilbert space $H_D$ with scalar product and norm $(u, v)_D = (D u, v), \ \| u \|_D = (u, v)_D^{1/2}$.

Let $A$ be a self-adjoint positive definite operator ($A: H \rightarrow H$):
\begin{equation}\label{1}
 A = A^*,
 \quad A \geq \delta I ,  
 \quad \delta > 0 , 
\end{equation} 
where $I$ is the identity operator in $H$.

We will consider the spectral problem
\[
 A \psi  = \lambda \psi  
\]
for the operator $A$. For the eigenvalues, we have
\[
 0 < \lambda_1 \leq \lambda_2 \leq \ldots \leq \lambda_K ,
\]
therefore in (\ref{1}) $\delta = \lambda_1$.
For each $u \in H$, we have the representation
\[
 u = \sum_{k=1}^{K} (u,\psi_k) \psi_k
\]
via eigenfunctions $\psi_k, \ \|\psi_k\| = 1, \ k = 1, 2, \ldots, K$.

We use the spectral definition of the operator function $f(A)$.
In this case, for the spectrum we have
\[
 \sigma(f(A)) = \{ f(\lambda): \ \lambda \in \sigma(A) \} ,
\] 
moreover
\[
 f(A) u = \sum_{k=1}^{K} f(\lambda_k) (u,\psi_k) \psi_k .
\] 
We consider the problem of finding the solution $u$ from the operator equation
\begin{equation}\label{2}
 f(A) u = \varphi .
\end{equation} 
For solvability, we will assume that $f(A)$ is a positive definite operator:
\begin{equation}\label{3}
 f(A) \geq \delta_f I, 
 \quad \delta_f = \min_{\lambda \in \sigma(A)} f(\lambda) > 0 . 
\end{equation} 
Special attention is paid to the problem of fractional powers of operators when
\begin{equation}\label{4}
 f(A) = A^\alpha ,
 \quad 0 < \alpha < 1 .  
\end{equation} 

We have well-developed numerical methods for solving equations with linear dependence of $f$ on $A$:
\[
 f(A) = c I + A,
 \quad c = \const,
 \quad c + \delta  > 0 . 
\] 
In particular, in many applied problems the operator $ A $ is associated with some approximation of a self-adjoint elliptic operator of the second order.
We want to build on this basis computational algorithms for nonlinear $f(A)$.

\section{Approximate solutions of equations with operator function} 

To solve problem (\ref{2}), we have
\begin{equation}\label{5}
 u = f^{-1}(A) \varphi .
\end{equation} 
Instead of $f^{-1}(A)$, some approximation 
\[
 r(A) \approx f^{-1}(A) 
\]
is used.
Thus, instead of $u$ we have an approximate solution
\begin{equation}\label{6}
 \accentset{\sim }{u} =  r(A) \varphi .
\end{equation}  

It is natural to associate the error estimates $\accentset{\sim }{u} - u$  with same estimates for the solution of the problem itself (\ref{2}). 
Taking into account (\ref{3}), for the problem (\ref{2}) it is possible, for example, to obtain a priori estimates in $H_D$, where $D = f (A), I$:
\begin{equation}\label{7}
 \|u\|_{f(A)} \leq \|\varphi \|_{f^{-1}(A)} \leq \delta_f^{-1/2}\|\varphi \|,
\end{equation}   
\begin{equation}\label{8}
 \|u\| \leq \|\varphi\|_{f^{-2}(A)} \leq \delta_f^{-1}\|\varphi \|.
\end{equation}   

For $\accentset{\sim }{u} - u$ we have
\begin{equation}\label{9}
 \accentset{\sim }{u} - u = \sum_{k=1}^{K} \big (r(\lambda_k) - f^{-1}(\lambda_k) \big )  (\varphi,\psi_k) \psi_k . 
\end{equation} 
We will assume that we have the absolute error of approximation of the function $f^{-1}(\lambda)$:
\begin{equation}\label{10}
 \big | r(\lambda) - f^{-1}(\lambda) \big | \leq \varepsilon ,
 \quad  \lambda \in [\lambda_1, \lambda_K ] .
\end{equation} 
For (\ref{10}) from (\ref{9}), we obtain
\begin{equation}\label{11}
  \|\accentset{\sim }{u} - u\|_D \leq \varepsilon \|\varphi \|_D ,
\end{equation} 
when $D = f(A), I$.
Comparison with (\ref{7}), (\ref{8}) shows that increased smoothness of the right-hand side is required.

When specifying the relative error
\begin{equation}\label{12}
 \big | r(\lambda) - f^{-1}(\lambda) \big | \leq \varepsilon f^{-1}(\lambda),
 \quad  \lambda \in [\lambda_1, \lambda_K ] ,
\end{equation} 
the situation is more favorable. In case (\ref{12}), we have
\begin{equation}\label{13}
  \|\accentset{\sim }{u} - u\|_{f(A)} \leq \varepsilon \|\varphi \|_{f^{-1}(A)} ,
\end{equation} 
\begin{equation}\label{14}
  \|\accentset{\sim }{u} - u\| \leq \varepsilon \|\varphi \|_{f^{-2}(A)} .
\end{equation} 

\subsection{Rational approximation} 

When constructing approximations for an operator function $f(A)$, we start from approximations of functions $f(\lambda)$.
Among the function approximations, we can single out the classical rational approximation \cite{braess1986nonlinear}.
Taking into account the computational implementation, we consider a special version of rational approximation, when the expression 
\begin{equation}\label{15} 
 r(A) = \sum_{i=1}^{m} a_i (b_i I + A)^{-1} 
\end{equation} 
is used for $r(A)$.

At this stage of the study, we will not somehow specify the choice of rational approximation (\ref{15}).
We will restrict ourselves only to the natural assumption about the positiveness of the coefficients:
\begin{equation}\label{16}
 a_i > 0, 
 \quad b_i > 0, 
 \quad i = 1, \ldots, m . 
\end{equation} 
Under these conditions, for the approximating operator  $r(A)$,  we have 
\begin{equation}\label{17}
 r(A) = r^*(A)  > 0 ,
 \quad r(A) A =  A r^*(A).  
\end{equation} 

Assumptions (\ref{17}) significantly narrows the class of problems (\ref{2}).
In particular, in this case the function $r(\lambda)$ is monotonically decreasing for $\lambda> 0$.
Therefore, we can count on the applicability of approximations (\ref{15}), (\ref{17}) for monotonically increasing functions $f(\lambda)$ for $\lambda> 0$.

For an approximate solution (\ref{6}) based on (\ref{15}), we have
\begin{equation}\label{18}
 \accentset{\sim }{u} =  \sum_{i=1}^{m} a_i w_i ,
\end{equation} 
where $w_i, \ i = 1, \ldots, m,$ are the solution to the problems
\begin{equation}\label{19}
 (b_i I + A) w_i = \varphi , 
 \quad i = 1, \ldots, m . 
\end{equation} 

We can slightly extend the class of rational approximations used by putting
\[
  r(A) = \sum_{i=1}^{m} a_i (b_i I + g(A))^{-1} ,
\] 
where
\begin{equation}\label{20}
 g(A) = (I + g_1 A) (g_2 I+ g_3A)^{-1} ,
 \quad g_n > 0, \quad n = 1,2,3.  
\end{equation} 
For $w_i, \ i = 1, \ldots, m,$ we get
\[
 (b_i (g_2 I+ g_3A)) + (I + g_1 A)) w_i = (g_2 I+ g_3A) \varphi , 
 \quad i = 1, \ldots, m .  
\] 

In general, problems (\ref{19}) are solved approximately.
It is necessary to adapt the accuracy of solving these subproblems to the peculiarities of our original problem (\ref{2}).
Let $\accentset{\sim }{w}_i,  \ i = 1, \ldots, m,$ are approximate solutions.
For an approximate solution of the problem (\ref{2}), we obtain
\begin{equation}\label{21}
 \accentset{\approx}{u} = \sum_{i=1}^{m} a_i \accentset{\sim }{w}_i ,
\end{equation} 
moreover
\begin{equation}\label{22}
 \|\accentset{\sim }{w}_i - w_i \| \leq \varepsilon_i \|\varphi \|,
 \quad i = 1, \ldots, m .  
\end{equation} 
We have
\[
 \|\accentset{\approx}{u} - \accentset{\sim }{u} \| \leq \sum_{i=1}^{m} a_i \varepsilon_i \|\varphi \| .
\] 
When specifying the error of the approximate solution for problems (\ref{19}) in the form
\begin{equation}\label{23}
 \varepsilon_i = \frac{\varepsilon_0}{b_i + \delta} ,
 \quad i = 1, \ldots, m ,
\end{equation} 
for the case (\ref{10}), we obtain
\[
\|\accentset{\approx}{u} - \accentset{\sim }{u} \| \leq \varepsilon_0 (f^{-1} (\delta) + \varepsilon) \|\varphi \| .
\] 
Taking into account (\ref{11}) ($D = I$), for the total error (approximation error $\accentset{\sim }{u} - u$ and error in solving problems (\ref{19})), we have the estimate
\begin{equation}\label{24}
\|\accentset{\approx}{u} - u \| \leq \big (\varepsilon + \varepsilon_0 (f^{-1} (\delta) + \varepsilon) \big ) \|\varphi \| .
\end{equation} 
This enables us to formulate the following statement.

\begin{thm}\label{t-1}
Let the estimate (\ref{10}) hold for the rational approximation (\ref{15}), (\ref{16}), and let the approximate solution $\accentset{\approx}{u}$ of the problem (\ref{5}) be determined according to (\ref{21}), where $\accentset{\sim }{w}_i,  \ i = 1, \ldots, m,$  are approximate solutions to problems (\ref{19}).
If the estimates (\ref{22}), (\ref{23}) are fulfilled, then the estimate (\ref{24}) holds for the error in solving the problem (\ref{5}).
\end{thm}

\subsection{Approximation by exponential sums} 

The second interesting class of function approximation \cite{braess1986nonlinear} is based on the sum of exponents: ES (Exponential Sums) approximation.
In this case, we have
\begin{equation}\label{25} 
 r(A) = \sum_{i=1}^{m} a_i \exp (- b_i A) .
\end{equation} 
We focus on approximations with positive coefficients (\ref{16}).
An approximate solution can be represented in the form (\ref{18}), where $w_i = w (b_i)$, and $w(t)$ is the solution to the Cauchy problem:
\begin{equation}\label{26}
 \frac{d w}{d t} + A w = 0,
 \quad 0 < t \leq T,
\end{equation} 
\begin{equation}\label{27}
 w(0) = \varphi ,
\end{equation} 
moreover
\[
 T = \max_{1 \leq i \leq m} b_i .
\] 

We can use a similar approach when applying the approximation
\begin{equation}\label{28}
 r(A) = \sum_{i=1}^{m} a_i \exp \big (- b_i g(A) \big ) ,  
\end{equation} 
where $g(A)$ is defined according to (\ref{20}).
In this case, instead of (\ref{26}), we use the equation 
\[
 (g_2 I + g_3 A) \frac{d w}{d t} + (I + g_1 A) w = 0,
 \quad 0 < t \leq T. 
\]  

Taking into account (\ref{1}), it is convenient to use instead of (\ref{25}) the representation
\[
 r(A) = \sum_{i=1}^{m} a_i \exp (- b_i \delta )  \exp \big (- b_i (A-\delta I) \big ) . 
\] 
Here, an approximate solution to the problem  (\ref{2}) is found from
\begin{equation}\label{29}
 \accentset{\sim }{u} =  \sum_{i=1}^{m} a_i \exp (- b_i \delta ) w_i ,
\end{equation} 
and to find $w_i = w (b_i)$ we solve the equation
\begin{equation}\label{30}
 \frac{d w}{d t} + (A - \delta I) w = 0,
 \quad 0 < t \leq T .
\end{equation} 

An approximate solution to the Cauchy problem (\ref{27}),  (\ref{30}) is denoted by $\accentset{\sim }{w}(t), \ 0 \leq t \leq T$.
We assume that there has been 
\begin{equation}\label{31}
 \|\accentset{\sim }{w}(t) - w(t) \| \leq \varepsilon_0 \|\varphi \|,
 \quad 0 < t \leq T .
\end{equation}  
Such estimates are discussed, for example, in the book \cite{thomee2006} using finite element approximations in space and finite difference approximations in time.
With considering
\begin{equation}\label{32}
 \accentset{\approx}{u} =  \sum_{i=1}^{m} a_i \exp (- b_i \delta ) \accentset{\sim }{w}_i , 
\end{equation} 
we will get
\[
\|\accentset{\approx}{u} - \accentset{\sim }{u} \| \leq \varepsilon_0 r(\delta) \|\varphi \| .
\] 
Under the assumption (\ref{10}), for the total error, we obtain the estimate (\ref{24}).

Similarly to Theorem \ref{t-1}, we can formulate the following statement about the error in the approximate solution of the problem (\ref{5}) when using the ES approximation.

\begin{thm}\label{t-2}
Let the estimate (\ref{10}) hold for the ES approximation (\ref{16}), (\ref{25}), and let the approximate solution $\accentset{\approx}{u}$ of the problem (\ref{5}) be determined according to (\ref{32}), where $\accentset{\sim }{w}_i = \accentset{\sim }{w}(b_i),  \ i = 1, \ldots, m,$  and $\accentset{\sim }{w}(t)$ is approximate solutions to problem (\ref{27}),  (\ref{30}).
If the estimate (\ref{31}) is fulfilled, then the estimate (\ref{24}) holds for the error in solving the problem (\ref{5}).
\end{thm}

\subsection{Approximation by exponential products} 

Let us note the possibilities that the method of approximating functions gives in the form of not a sum of exponents but a product of exponentials: EP (Exponential Products) approximation.
As we see it, this is a new class of function approximation that has not been considered previously.

We will start from the record
\begin{equation}\label{33}
 f^{-1}(A) = \exp \big (- \log \big (f(A) \big ) \big )  
\end{equation} 
and approximations of the function $\log \big (f(A) \big )$.
Similar to rational approximation (\ref{15}) for $f^{-1}(A)$, we put
\begin{equation}\label{34}
 \log \big (f(A) \big ) \approx s(A) = a_0 I + \sum_{i=1}^{m} (a_i I + b_i A) (c_i I + d_i A)^{-1} .
\end{equation} 
We will assume that for the coefficients $a_i, b_i, c_i, d_i, \ i = 1,\dots, m,$ constraints of the type (\ref{16}) are again satisfied. Taking into account (\ref{1}), we formulate assumptions in the form 
\begin{equation}\label{35}
 a_i - \delta b_i \geq 0, 
 \quad c_i - \delta d_i > 0, 
 \quad i = 1, \ldots, m . 
\end{equation} 
Thus, we apply the EP approximation
\begin{equation}\label{36}
 r(A) = \exp \Big (- \sum_{i=1}^{m} (a_i I + b_i A) (c_i I + d_i A)^{-1} \Big ) \exp (-a_0 I) . 
\end{equation}  

For an approximate solution (\ref{6}), (\ref{36}), we obtain
\begin{equation}\label{37}
 \accentset{\sim }{u} = w_m(1) ,
\end{equation} 
when solving a sequence of equations:
\begin{equation}\label{38}
 (c_i I + d_i A) \frac{d w_i}{d t} + (a_i I + b_i A) w_i = 0,
 \quad 0 < t \leq 1 .
\end{equation} 
For (\ref{35}), we have
\[
 a_i I + b_i A \geq  0 ,
 \quad c_i I + d_i A > 0,
 \quad i = 1, \ldots, m . 
\] 
The initial conditions for (\ref{38}) are of the form
\begin{equation}\label{39}
 w_1(0) =  \exp (-a_0) \varphi ,
 \quad w_i(0) = w_{i-1} (1),
 \quad i = 2,\ldots, m .  
\end{equation} 

Instead of sequentially solving problems (\ref{38}), (\ref{39}), you can use one equation with piecewise constant coefficients.
Let $w(t)$ satisfy the equation
\begin{equation}\label{40}
 \big (c(t) I + d(t) A \big ) \frac{d w}{d t} + \big (a(t) I + b(t) A \big ) w = 0,
 \quad 0 < t \leq m ,
\end{equation} 
in which
\[
 \big (a(t), b(t), c(t), d(t) \big ) = \big (a_i, b_i, c_i, d_i \big ),
 \quad i-1 < t \leq i,
 \quad i = 1,\ldots, m .   
\] 
The equation (\ref{40}) is supplemented with the initial condition 
\begin{equation}\label{41}
 w(0) =  \exp (-a_0) \varphi .
\end{equation} 
Taking into account (\ref{37}), the approximate solution of the problem (\ref{6}), (\ref{35}) is
\begin{equation}\label{42}
  \accentset{\sim }{u} = w(m) .
\end{equation} 

The main difference between EP approximation and ES approximation is related to the approximation of not the operator function $f^{-1}(A)$ (see (\ref{25})), but the function $\log \big (f(A) \big )$ (see (\ref{34})).
An approximate solution is constructed on the basis of the Cauchy problem for a slightly more complicated equation (compare (\ref{30}) and (\ref{40})).

Let us define the approximation error $\log \big (f(A) \big )$ similarly to (\ref{10}):  
\begin{equation}\label{43}
 \big | s(\lambda ) - \log \big (f(\lambda ) \big ) \big | \leq \varepsilon_1 ,
 \quad  \lambda \in [\lambda_1, \lambda_K ] . 
\end{equation} 
Taking into account (\ref{36}), we have
\[
 r(\lambda) - f^{-1}(\lambda) = \Big (\exp\Big ( - \Big (s(\lambda) - \log \big (f(\lambda ) \big ) \Big ) - 1  \Big ) f^{-1}(\lambda) ,
 \quad  \lambda \in [\lambda_1, \lambda_K ] .
\] 
For (\ref{43}) and for small $\varepsilon_1$, we obtain
\begin{equation}\label{44}
 \big | r(\lambda) - f^{-1}(\lambda) \big | \leq \big (\exp(\varepsilon_1) - 1 \big ) f^{-1}(\lambda),
 \quad  \lambda \in [\lambda_1, \lambda_K ] .
\end{equation} 
Thus, we obtain an estimate of the relative error (\ref{12}) when
\begin{equation}\label{45}
 \varepsilon = \exp(\varepsilon_1) - 1 \simeq \varepsilon_1 . 
\end{equation} 
We have an estimate of the relative accuracy of the approximation.
In view of this, we can use the estimates of the accuracy (\ref{13}), (\ref{14}) for the approximate solution of the problem (\ref{5}).

Investigation of the influence of the error in solving the problem (\ref{40}), (\ref{41}) requires a separate consideration.
We denote the approximate solution of the problem (\ref{40}), (\ref{41}) by $\accentset{\sim }{w}(t), \ 0 < t \leq m$.
For an approximate solution of the problem (\ref{5}), from (\ref{42}), we obtain
\begin{equation}\label{46}
  \accentset{\approx}{u} = \accentset{\sim }{w}(m) .
\end{equation} 
We will assume that, like (\ref{31}), the estimate 
\begin{equation}\label{47}
 \|\accentset{\sim }{w}(m) - w(m) \| \leq \varepsilon_0 \|\varphi \| 
\end{equation}   
holds.
For $\accentset{\approx}{u} - u$, from (\ref{14}) and (\ref{46}), (\ref{47}), we obtain the estimate
\begin{equation}\label{48}
 \|\accentset{\approx}{u} - u \| \leq \varepsilon_0  \|\varphi \| 
 + \big (\exp(\varepsilon_1) - 1 \big )  \|\varphi \|_{f^{-2}(A)} .
\end{equation} 

\begin{thm}\label{t-3}
Let the estimate (\ref{43}) hold for approximation errors $\log \big (f(\lambda )$, and let the approximate solution $\accentset{\approx}{u}$ of the problem (\ref{5}) be determined according to (\ref{46}), where $\accentset{\sim }{w}(t), \ 0 < t \leq m,$ is approximate solutions to problem (\ref{40}), (\ref{41}). 
If the estimate (\ref{47}) is fulfilled, then the estimate (\ref{48}) holds for the error in solving the problem (\ref{5}).
\end{thm}

\section{Approximations of the fractional power of an operator based on the integral representation} 

We illustrate the approaches to approximating the operator function described by us on a typical problem.
As such, we will take the case (\ref{4}); problems with a fractional power of an elliptic operator are actively discussed in the literature.
The main methods for the numerical solution of such problems are investigated, for example, in \cite{harizanov2020rev,bonito2018numerical}.
In the case (\ref{4}), for the solution of the equation (\ref{2}), we have
\begin{equation}\label{49}
 u = A^{-\alpha} \varphi ,
\end{equation} 
for $0 < \alpha < 1$.

\subsection{Rational approximation}  

Usually (see, for example, \cite{lunardi2018interpolation,yagi2009abstract}) for a fractional power of a positive definite self-adjoint operator, one uses the integral representation
\begin{equation}\label{50}
 A^{-\alpha} = \frac{\sin(\alpha \pi)}{\pi} \int_{0}^{\infty } \theta^{-\alpha} (\theta I + A )^{-1} d \theta .
\end{equation}
This Balakrishnan formula \cite{balakrishnan1960fractional} follows from the Dunford integral representation \cite{yoshida1995functional} of the operator function when
\[
 f(A) = \frac{1}{2 \pi i} \int_{\partial U} f(\lambda) (\lambda I - A)^{-1} d \lambda ,
 \quad \sigma (A) \subseteqq U . 
\]  

Approximation of $A^{-\alpha}$  is carried out on the basis of using some quadrature formula for the integral on the right-hand side (\ref{50}).
In this case,
\begin{equation}\label{51}
 A^{-\alpha} \approx  \frac{\sin(\alpha \pi)}{\pi}  \sum_{i=1}^{m} w_i \theta_i^{-\alpha} (\theta_i I + A)^{-1} ,
\end{equation}
where $w_i$ are the weights of the quadrature formula, and $\theta_i, \ i =1, \ldots, m,$ are the nodes of the quadrature formula.
Taking into account (\ref{51}) for the coefficients of rational approximation (\ref{15}), we obtain
\[
 a_i =  \frac{\sin(\alpha \pi)}{\pi} w_i \theta_i^{-\alpha} ,
 \quad b_i =  \theta_i,
 \quad  i =1, \ldots, m .  
\] 

The problems of constructing quadratures are caused, on the one hand, by the peculiarity of the integrand at $\theta \rightarrow 0$ and the unboundedness of the integration interval, on the other hand.
To simplify the problem of constructing quadrature formulas, new variables are used.
For example, to eliminate the singularity at $\theta \rightarrow 0$, we can
(for example \cite{bonito2015numerical}) use the new variable $\eta$: $\theta = \exp(\eta)$.
Thus, from (\ref{50}), we obtain
\[
 A^{-\alpha} = \frac{\sin(\alpha \pi)}{\pi} \int_{-\infty }^{\infty } \exp \big ((1-\alpha)\eta\big) \big (\exp(\eta) I + A \big)^{-1} d \eta .
\]
The peculiarities of the behavior of the integrand in this case are taken into account by choosing sinc numerical methods for quadrature \cite{stenger2011handbook}. 

To obtain the integral on a finite interval, in \cite{frommer2014efficient,Aceto2019}, there is employed the transformation
\[
 \theta = \mu \frac{1-\eta}{1+\eta},
 \quad \mu > 0.
\] 
From (\ref{50}), we have
\[
  A^{-\alpha  } = \frac{2 \mu^{1-\alpha } \sin(\pi \alpha  )}{\pi} \int_{-1}^{1} (1-\eta)^{-\alpha }(1+\eta)^{\alpha -1}
  \big (\mu (1-\eta) I + (1+\eta) A \big )^{-1} d \eta . 
\]
To approximate the right-hand side, we apply the Gauss-Jacobi quadrature formula \cite{Rabinowitz} 
with the weight $(1-\eta)^{-\alpha }(1+\eta)^{\alpha -1}$.

Using
\[
 \theta = \frac{\eta^{1/(1-\alpha)}}{(1-\eta)^{\varkappa/\alpha}} ,
\] 
the formula (\ref{50}) takes the form
\begin{equation}\label{52}
  A^{-\alpha } = \frac{\sin(\pi \alpha )}{(1-\alpha)\pi} \int_{0}^{1} (1-\eta)^{\varkappa -1} 
    \Big (1+ \Big(\varkappa\frac{1-\alpha}{\alpha} - 1 \Big ) \eta \Big )
    \Big (\eta^{\frac{1}{1-\alpha}} I + (1-\eta)^{\frac{\varkappa}{\alpha}} A \Big )^{-1} d \eta   
\end{equation} 
with the parameter $\varkappa > 1$. 
We used such an integral representation in the article \cite{vabishchevich2020}.
The advantages of this representation are (i) we integrate on a finite interval, (ii) we avoid the singularity of the integrand, (iii) and we control the smoothness of the integrand by choosing the parameter $\varkappa$ of the integral representation.
Because of this, we can restrict ourselves to standard quadrature formulas.

\subsection{ES approximation} 

Boundary value problems for fractional powers of elliptic operators are often considered \cite{stinga2019user} in terms of methods of semigroups.
In this case, the integral representation of the solution to the problem (\ref{49}) takes the form
\begin{equation}\label{53}
  A^{-\alpha} = \frac{1}{\Gamma(\alpha)} 
	\int_0^\infty \theta^{\alpha-1} \exp(- \theta  A) \, d \theta ,
\end{equation}
where $\Gamma(\alpha)$ is the gamma function. 

When using a suitable quadrature formulas we obtain from (\ref{53}) to
\begin{equation}\label{54}
 A^{-\alpha} \approx \frac{1}{\Gamma(\alpha)} \sum_{i=1}^{m} w_i  \exp(- \theta_i  A) .
\end{equation}
Taking into account (\ref{54}), the ES approximation coefficients (\ref{25}) are expressed in terms of the weights and nodes of the quadrature formula as follows
\[
 a_i =  \frac{1}{\Gamma(\alpha)} w_i ,
 \quad b_i =  \theta_i,
 \quad  i =1, \ldots, m .  
\] 
The works \cite{cusimano2018discretizations,cusimano2020numerical} apply the quadrature formula weights are associated with the function $\theta^{\alpha-1}$; the authors distinguish regions of integration near $\eta = 0$ and for large $\eta$.

For a self-adjoint positive definite operator $A$ (see (\ref{1})), instead of (\ref{53}), you can use the representation
\begin{equation}\label{55}
  A^{-\alpha} = \frac{1}{\Gamma(\alpha)} 
	\int_0^\infty \theta^{\alpha-1} \exp(-\delta \theta)  \exp\big (- \theta  (A-\delta I) \big ) \, d \theta .
\end{equation}
The peculiarity of the integrand $\theta^{\alpha-1} \exp(-\delta \theta)$ in (\ref{55}) is most fully taken into account using the generalized Gauss-Laguerre quadrature formula \cite{davis2007methods}.
In this case, we have
\[
 A^{-\alpha} \approx \frac{1}{\Gamma(\alpha)} \sum_{i=1}^{m} w_i  \exp \big (- \theta_i  (A-\delta I) \big ) .
\]
We get the approximation (\ref{28}) with
\[
 a_i =  \frac{1}{\Gamma(\alpha)}\exp(\delta \theta_i) \, w_i ,
 \quad b_i =  \theta_i,
 \quad  i =1, \ldots, m .  
\] 
Such ES approximations for the fractional power of the operator are constructed in our article \cite{vabishchevich2021fractional}.

\subsection{EP approximation} 

In the case of a fractional power of the operator $A$, we have
\[
 \log(A^{\alpha}) = \alpha \log(A) .
\] 
For $\lambda \geq \delta$, we obtain
\begin{equation}\label{56}
 \log(\lambda) = \log(\delta) + \int_{0}^{1} \frac{\lambda - \delta}{\theta (\lambda - \delta) + \delta } d \theta .
\end{equation} 
Taking into account (\ref{56}) for $\log(A^{\alpha})$, we have the representation
\begin{equation}\label{57}
 \log(A^{\alpha}) = \alpha \log(\delta) I + \alpha \int_{0}^{1} (A - \delta I ) \big (\theta (A - \delta I ) + \delta I \big )^{-1} d \theta .
\end{equation} 
A similar representation for matrix logarithm ($\alpha = 1, \ \delta = 1$ in (\ref{57})) is introduced by Richter \cite{richter1949logarithmus} (see also \cite{higham2008functions,wouk1965integral}).

Using some quadrature formula, we obtain
\[
  \log(A^{\alpha}) \approx \alpha \log(\delta) I + \alpha  \sum_{i=1}^{m} w_i (A - \delta I ) \big (\theta_i (A - \delta I ) + \delta I \big )^{-1} .
\] 
For the representation coefficients (\ref{34}), we have
\[
\begin{split}
 a_0 & = \alpha \log(\delta),
 \quad a_i = - \alpha w_i \delta ,
 \quad b_i = \alpha w_i , \\
 c_i & = \delta (1 - \theta_i),
\quad  d_i = \theta_i, 
\quad  i = 1, \ldots, m . 
\end{split}
\] 

When constructing operator approximations based on the Balakrishnan formula (\ref{50}), great attention is paid to the use of new integration variables. With this approach, we can expect to improve accuracy when using both the ES approximation and the EP approximation.

\section{Cauchy method for problems with a fractional power of an operator} 

In this part of the work, we discuss the possibilities of an approximate solution of problems with a fractional power of an operator based on the transition to an auxiliary evolutionary problem.
The use of some approximations for the numerical solution of the corresponding Cauchy problem can be interpreted as some rational approximations.
The starting point for us is the integral representation for the operator function.
The upper limit of integration is variable and differentiation with respect to this parameter gives the corresponding equation for the desired solution.

\subsection{Rational approximation} 

According to (\ref{50}), we will consider the function
\begin{equation}\label{58}
 v(t) = \frac{\sin(\alpha \pi)}{\pi} \int_{0}^{t} \theta^{-\alpha} (\theta I + A )^{-1} d \theta \varphi .
\end{equation} 
Then the solution to problem (\ref{49}) is $u = v(\infty)$. 
Differentiating (\ref{58}) with respect to $t$, we obtain
\begin{equation}\label{59}
 t^{\alpha} (t I + A ) \frac{d v}{d t} = \frac{\sin(\alpha \pi)}{\pi}  \varphi .
\end{equation} 
For this ordinary differential equation, the initial condition is used
\begin{equation}\label{60}
 v(0) = 0 .
\end{equation} 

The singularity in the equation (\ref{59}) can be removed by replacing $t \rightarrow t^{1+\alpha}$. 
When using the new variables of integration, we have for $A^{-\alpha}$ the integral representation
\begin{equation}\label{61}
 A^{-\alpha} = \frac{\sin(\alpha \pi)}{\pi} \int_{0}^{\infty } (\theta^{\frac{1}{1+\alpha} }  I + A )^{-1} d \theta .
\end{equation} 
In this case, for $v(t)$, we use the equation
\begin{equation}\label{62}
  \big (t^{\frac{1}{1+\alpha} } I + A \big ) \frac{d v}{d t} = \frac{\sin(\alpha \pi)}{\pi}  \varphi .
\end{equation}
Integral representation (\ref{52}) leads us to
\[
  v(t) = \frac{\sin(\pi \alpha )}{(1-\alpha)\pi} \int_{0}^{t} (1-\eta)^{\varkappa -1} 
    \Big (1+ \Big(\varkappa\frac{1-\alpha}{\alpha} - 1 \Big ) \eta \Big )
    \Big (\eta^{\frac{1}{1-\alpha}} I + (1-\eta)^{\frac{\varkappa}{\alpha}} A \Big )^{-1} d \eta \, \varphi .  
\] 
This gives the evolutionary equation
\[
  \Big (t^{\frac{1}{1-\alpha}} I + (1-t)^{\frac{\varkappa}{\alpha}} A \Big ) \frac{d v}{d t} = 
 (1-t)^{\varkappa -1} \Big (1+ \Big(\varkappa\frac{1-\alpha}{\alpha} - 1 \Big ) t \Big )\frac{\sin(\pi \alpha )}{(1-\alpha)\pi} \varphi . 
\] 
The coefficients of the equation have no singularities and the Cauchy problem with the condition (\ref{60}) is solved on the unit interval, with  $u = v(1)$.

The use of certain approximations in time leads us to some rational approximation.
Let, for example, $0 = t_0 < t_1, ..., $ are grid nodes in time.
For an approximate solution of the problem (\ref{60}), (\ref{62}), you can use the simplest two-level approximations, when
approximate solution $v_n \approx v(t_n)$  is determined from
\begin{equation}\label{63}
 \Big (t_{n+1/2}^{\frac{1}{1+\alpha} } I + A \Big ) \frac{v_{n+1}-v_n}{t_{n+1}-t_n} = \frac{\sin(\alpha \pi)}{\pi}  \varphi ,
 \quad n = 0, 1, \ldots ,  
\end{equation} 
where, for example, $t_{n+1/2} = (t_{n+1}+t_n)/2$ for $v_0 = 0$. 
The difference scheme (\ref{63}) corresponds to the use of rational approximation with the coefficients
\[
 a_i =  \frac{\sin(\alpha \pi)}{\pi} (t_{i+1}-t_i) ,
 \quad b_i =  t_{i+1/2}^{\frac{1}{1+\alpha} } ,
 \quad  i =1, \ldots, m .  
\] 
The approximation (\ref{63}) is associated with the simplest quadrature rectangle formula for the integral representation (\ref{61}).

\subsection{ES approximation} 

Taking into account (\ref{53}), we introduce the function  $v(t)$ in the form
\begin{equation}\label{64}
  v(t) = \frac{1}{\Gamma(\alpha)} 
	\int_0^t \theta^{\alpha-1} \exp(- \theta  A) \, d \theta \, \varphi ,
\end{equation}
so $u = v(\infty)$.
Differentiation (\ref{64}) with respect to $t$ gives the system of equations
\begin{equation}\label{65}
 \Gamma(\alpha) t^{1-\alpha} \frac{d v}{d t} - w = 0, 
\end{equation} 
\begin{equation}\label{66}
 \frac{d w}{d t} + A w = 0 ,
 \quad t > 0 .  
\end{equation} 
This system of equations is supplemented by the initial conditions
\begin{equation}\label{67}
 v(0) = 0,
 \quad w(0) = \varphi  . 
\end{equation} 

As we noted above (see (\ref{29}), (\ref{30})), it is appropriate to go from $A$ to $A - \delta I$ when
\[
  v(t) = \frac{1}{\Gamma(\alpha)} 
	\int_0^t \theta^{\alpha-1} \exp(-\delta \theta) \exp \big (- \theta  (A - \delta I) \big ) \, d \theta \, \varphi . 
\] 
In this case, the system of equations (\ref{65}), (\ref{66}) takes the form
\[
 \Gamma(\alpha) t^{1-\alpha} \frac{d v}{d t} - \exp(-\delta t) w = 0, 
\]
\[
 \frac{d w}{d t} + (A - \delta I) w = 0 ,
 \quad t > 0 .  
\]

The problem (\ref{65})--(\ref{67}) can be solved sequentially: first $w(t)$ and then $v(t)$ using some approximations in $t$.
ES approximation is associated with partial discretization when $t$ approximation is performed for only one equation (\ref{65}).
For example, the simplest approximation
\[
 \Gamma(\alpha) t_{n+1/2}^{1-\alpha} \frac{v_{n+1}-v_n}{t_{n+1}-t_n} = w(t_{n+1/2}),
 \quad n = 0, 1, \ldots ,  
\] 
corresponds to ES approximation (\ref{53}) when specifying coefficients
\[
 a_i =  \frac{1}{\Gamma(\alpha)} t_{i+1/2}^{\alpha-1} (t_{i+1}-t_i) ,
 \quad b_i =  t_{i+1/2} ,
 \quad  i =1, \ldots, m .  
\] 
In this case, the quadrature rectangle formula for the integral representation (\ref{52}) is applied.

Let us note some other possibilities of an approximate solution based on the transition to the Cauchy problem.
Taking into account the integral representation (\ref{52}) using the change of variables $\theta \rightarrow \theta^\alpha$, we obtain
\[
  A^{-\alpha} = \frac{1}{\alpha \Gamma(\alpha)} 
	\int_0^\infty \exp(- \theta^{\frac{1}{\alpha}} A) \, d \theta .
\]
We now define the function $v(t)$ by the expression
\begin{equation}\label{68}
 v(t) = \frac{1}{\alpha \Gamma(\alpha)} 
	\int_0^t \exp \big (- \theta^{\frac{1}{\alpha}} A \big ) \, d \theta \, \varphi .
\end{equation} 
For the first derivative we have
\begin{equation}\label{69}
 \frac{d v}{d t} = \frac{1}{\alpha \Gamma(\alpha)} 
	\exp \big (- t^{\frac{1}{\alpha}} A \big ) \, \varphi .
\end{equation}  
Repeated differentiation leads us to the second-order equation
\begin{equation}\label{70}
 \frac{d^2 v}{d t^2}  + t^{\frac{1-\alpha }{\alpha}} A \frac{d v}{d t} = 0 .
\end{equation}  
From (\ref{68}), (\ref{69}) we have the initial conditions for the equation (\ref{70}) 
\begin{equation}\label{71}
 v(0) = 0,
 \quad \frac{d v}{d t}(0) =  \frac{1}{\alpha \Gamma(\alpha)} \varphi .
\end{equation} 
A solution to the Cauchy problem (\ref{70}), (\ref{71}) for $t \rightarrow \infty$ gives a solution to the problem (\ref{49}).

\subsection{EP approximation} 

Taking into account the integral representation (\ref{57}), we define the function
\begin{equation}\label{72}
 v(t) = \exp \Big (- \alpha \log(\delta) I - 
 \alpha \int_{0}^{t} (A - \delta I ) \big (\theta (A - \delta I ) + \delta I \big )^{-1} d \theta \Big ) \varphi . 
\end{equation} 
The solution to problem (\ref{49}) is $u = v(1)$.
Differentiation (\ref{72}) with respect to $t$ leads us to the equation
\begin{equation}\label{73}
 \big (\theta (A - \delta I ) + \delta I \big )\frac{d v}{d t} + \alpha (A - \delta I ) v = 0 .
\end{equation} 
Taking into account (\ref{72}), the initial condition for (\ref{73}) is
\begin{equation}\label{74}
 v(0) = \delta^{-\alpha} \varphi .
\end{equation} 
With an approximate solution, the equation (\ref{73}) can be associated with some version of the equation (\ref{38}); this version arises when approximated by exponential products.

The Cauchy problem method based on the approximate solution of the problem (\ref{73}), (\ref{74}) was used in \cite{vabishchevich2014numerical}  for problems with a fractional power of an elliptic operator. 
To improve the accuracy of approximate solution, two- and three-level schemes of higher-order approximation on regular and irregular
grids are used (see e.g. Refs. \cite{duan2018numerical,CegVab2019,vciegis2020high}).

\section{Conclusion} 

\begin{enumerate} 
\item General approaches to the approximate solution of equations with functions of self-adjoint positive definite operators are considered.
The methods are based on the use of various types of function approximations and the solution of auxiliary problems.
The classical rational approximations and approximations by sums of exponentials are noted.
 \item A new class of function approximations by the product of exponentials is proposed. In the approximate solution of equations with operator functions, the computational implementation is based on the solution of first-order evolution equations.
 \item Estimates are obtained for the proximity of an approximate solution of equations with operator functions for rational approximation, approximation by sums of exponentials or their products. Variants with absolute and relative errors of approximation of operator functions are highlighted separately.
 \item General approaches are illustrated by constructing an approximate solution to problems with a fractional power of an operator.
Rational approximation, approximation by sums or by the product of exponentials, is based on the use of some integral representations of operator functions.
 \item Various approximations are used to construct an approximate solution to a problem with a fractional power of an operator based on the auxiliary Cauchy problem for some evolutionary equation.
\end{enumerate} 

\section*{Acknowledgements}

The publication has been prepared with the support by the mega-grant of the Russian Federation Government 14.Y26.31.0013
and the research grant 20-01-00207 of Russian Foundation for Basic Research.

\end{document}